\documentclass[12pt,twoside]{amsart}
\usepackage{amssymb}
\usepackage{graphicx}

\usepackage[margin=1in]{geometry}
\numberwithin{equation}{section}
\newcommand{\De}{\Delta}
\newcommand{\Ga}{\Gamma}
\newcommand{\si}{\sigma}
\newcommand{\Si}{\Sigma}
\newcommand{\NN}{\mathbb{N}}
\newcommand{\ZZ}{\mathbb{Z}}
\newcommand{\dup}{\overline{d}}
\newcommand{\ddown}{\underline{d}}

\DeclareMathOperator{\link}{link}
\DeclareMathOperator{\fdel}{del}

\newcommand{\st}{\; | \;}
\numberwithin{figure}{section}

\newtheorem{theorem}{Theorem}[section]
\newtheorem{lemma}[theorem]{Lemma}
\newtheorem{proposition}[theorem]{Proposition}
\newtheorem{corollary}[theorem]{Corollary}

\theoremstyle{definition}
\newtheorem{definition}[theorem]{Definition}
\newtheorem{remark}[theorem]{Remark}
\newtheorem{example}[theorem]{Example}

\newtheorem*{acknowledgement}{Acknowledgement}

\begin{document}

\title{The structure of the Boij-S\"oderberg posets}
\author[D.\ Cook II]{David Cook II}
\address{Department of Mathematics, University of Kentucky, 715 Patterson Office Tower, Lexington, KY 40506-0027, USA}
\email{dcook@ms.uky.edu}
\thanks{Part of the work for this paper was done while the author was partially supported by the National Security Agency under Grant Number H98230-09-1-0032.}

\begin{abstract}
    Boij and S\"oderberg made a pair of conjectures, which were subsequently proven by Eisenbud and Schreyer and then extended by Boij and 
    S\"oderberg, about the structure of Betti diagrams of Graded modules.  In the theory, a particular family of posets, and their associated
    order complexes, play an integral role.  We explore the structure of this family.  In particular, we show the posets are bounded complete
    lattices and the order complexes are vertex-decomposable, hence Cohen-Macaulay and squarefree glicci.
\end{abstract}

\maketitle

Boij and S\"oderberg recently conjectured in~\cite{BS08a} a complete characterisation, up to multiplication by a positive rational, of the structure of
Betti diagrams of finitely generated graded modules.  Their conjectures were proven for the Cohen-Macaulay case in~\cite{ES}.  These were further extended
to arbitrary graded modules in~\cite{BS08b}.  This characterisation centered around constructing the convex hull of a particular class of Betti diagrams and
showing it is equal to a geometric realisation of a particular simplicial complex--the order complex of a poset. 

The family of posets described in~\cite[Definition~2.3]{BS08a}, which we refer to as the {\em Boij-S\"oderberg posets}, is the focus of this paper.
In particular, we study the structure of the posets and their associated order complexes.  While no immediate applications are presented here, we
hope that our results will help shed light on the mysterious relation between the decomposition of Betti tables and the associated modules.

In Section~\ref{sec:preliminaries} we recall the relevant combinatorial definitions and explicitly define the Boij-S\"oderberg posets
(Definition~\ref{def:boij-soederberg-posets}).  Following this, in Section~\ref{sec:basic-structure} we determine some basic structural results for
the posets.  In particular, we show that the Boij-S\"oderberg posets are indeed bounded complete lattices (Proposition~\ref{pro:lattice}).
And in Section~\ref{sec:recursive-structure} we discuss the recursive structure of the posets by finding a recursive atom ordering
(see~\cite{BW83}) for each of the posets (Theorem~\ref{thm:atom-ordering}). This allows us to conclude that the order complexes are
vertex-decomposable, Cohen-Macaulay, and squarefree glicci (Corollary~\ref{cor:vd}).

\section{Preliminaries} \label{sec:preliminaries}

\subsection{Combinatorics}

A {\em simplicial complex} $\De$, on a finite set $V$, is a set of subsets of $V$ closed under inclusion; elements of $\De$ are called
{\em faces}.  The {\em dimension} of a face $\si$ is $\#\si - 1$ and of a complex $\De$ is the maximum of the dimensions of its faces.  A
complex whose maximal faces, called {\em facets}, are equi-dimensional is called {\em pure} and a complex with a unique maximal face is
called a {\em simplex}.

Given two simplicial complexes $\De$ and $\Ga$ with disjoint vertex sets, we define the {\em join} of $\De$ and $\Ga$ to be the
simplicial complex 
\[
    \De \star \Ga := \{ \si \cup \tau \st \si \in \De, \tau \in \Ga \}.
\]
If $\De = v$ is a vertex, then $v \star \Ga$ is said to be the {\em cone} of $\Ga$ with {\em apex} $v$.

Let $\si$ be a face of $\De$, then the {\em link} and {\em deletion} of $\si$ from $\De$ are given by
\[
    \link_\De{\si} := \{ \tau \in \De \st \tau \cap \si = \emptyset, \tau \cup \si \in \De \}
    \mbox{ and }
    \fdel_\De{\si} := \{ \tau \in \De \st \si \nsubseteq \tau \}.
\]
Following~\cite[Definition~2.1]{PB}, a pure complex $\De$ is said to be to be {\em vertex-decomposable} if either $\De$ is a
simplex or there exists a vertex $v \in \De$, called a {\em shedding vertex}, such that both $\link_{\De}{v}$ and $\fdel_{\De}{v}$ are
vertex-decomposable.  Checking if a particular simplicial complex is vertex-decomposable can be done using a computer program such
as~\cite{M2}, in particular, the package described in~\cite{Co} provides the appropriate methods.

A {\em poset} $P$ is a set with a partial ordering, that is, a binary relation ``$\leq$'' over the set which is reflexive, antisymmetric,
and transitive.  Given $v,u \in P$ with $v \leq u$, the {\em interval} $[v,u]$ is the sub-poset $\{w \in P \st v \leq w \leq u]$.
A {\em chain} is a sequence $v_0 < v_1 < \cdots < v_p$ of elements of $P$; such a chain is said to have {\em length} $p$.
The poset $P$ is called {\em pure} if every maximal chain has the same length.  The poset $P$ is {\em bounded} if there exists a unique minimal
element, $\hat{0}$, and a unique maximal element, $\hat{1}$, of $P$.

Let $S$ be a subset of $P$.  The {\em meet} of $S$, if it exists, is the infimum of $S$ and is denoted $\vee S$.  Similarly, the {\em join} of $S$,
if it exists, is the supremum of $S$ and is denoted $\wedge S$.  If every distinct pair of elements in $P$ has a meet and a join, then $P$ 
is called a {\em lattice}; if every subset of $P$ has a meet and a join, then $P$ is called a {\em complete} lattice.

Let $P$ be a finite poset.  For elements $x, y \in P$, we say that $y$ {\em covers} $x$, denoted $x \rightarrow y$, if $x < y$ and
$x < z \leq y$ implies $y = z$; in this case, we also say $x$ is {\em covered by} $y$.  An {\em atom} of a bounded poset is an element
which covers $\hat{0}$, the unique minimal element of $P$.

In~\cite{BW83}, Bj\"orner and Wachs define a bounded pure poset $P$ to {\em admit a recursive atom ordering} if either $P$ has
maximal chains of length one or there is an ordering of the atoms, $a_1, \ldots, a_t$, of $P$ which satisfies the properties:
\begin{enumerate}
    \item for $1 \leq j \leq t$, $[a_j, \hat{1}]$ admits a recursive atom ordering on its atoms $b_1, \ldots, b_s$ with the property that
        there exists a $1 \leq k \leq s$ such that $b_l \rightarrow a_i$ for some $1 \leq i < j$ if and only if $l \leq k$, and
    \item for $1 \leq i < j \leq t$, if $a_i, a_j < y$, then there is a $k < j$ and $z \in P$ such that $a_k, a_j \rightarrow z \leq y$.
\end{enumerate}
      
Associated to every finite poset $P$ is the {\em order complex}, denoted $\De(P)$, which is a simplicial complex with faces given by
chains in $P$; the facets of $\De(P)$ are exactly the maximal chains of $P$, hence $P$ is pure if and only if $\De(P)$ is pure.  Notice
that the minimal non-faces of $\De(P)$ are exactly the pairs of incomparable elements of $P$, so $\De(P)$ is a {\em flag complex}.

\subsection{Boij-S\"oderberg posets and order complexes}

In~\cite{BS08a}, Boij and S\"oderberg made a pair of conjectures about the possible graded Betti numbers of graded modules
up to multiplication by positive rational numbers; the conjecture was proven in~\cite{ES} and~\cite{BS08b}.  In the
course of their construction, they define a family of posets which we recall here.

\begin{definition} \label{def:boij-soederberg-posets}
    Let $\ddown, \dup \in \ZZ^{p+1}$ be strictly increasing sequences with $\ddown_i \leq \dup_i$ for $0 \leq i \leq p$.  Define
    the {\em Boij-S\"oderberg poset} of $\ddown$ and $\dup$ to be the set $\Pi_{\ddown, \dup}$ of strictly increasing sequences
    $d \in \ZZ^{p+1}$ such that $\ddown_i \leq d_i \leq \dup_i$ for $0 \leq i \leq p$ endowed with the partial ordering defined by
    $d \leq e$ whenever $d_i \leq e_i$ for $0 \leq i \leq p$.

    The {\em Boij-S\"oderberg order complex} of $\ddown$ and $\dup$ is the order complex of the Boij-S\"oderberg poset $\Pi_{\ddown, \dup}$,
    that is, $\De(\Pi_{\ddown, \dup})$.
\end{definition}

Notice that in~\cite{BS08a}, the poset $\Pi_{\ddown, \dup}$ has the associated pure diagrams as vertices, but these are equivalent as
pure diagrams are in bijection to strictly increasing sequences in $\ZZ^{p+1}$.  We further note that $\De(\Pi_{\ddown, \dup})$ is
a pure complex by~\cite[Proposition~2.7]{BS08a}.

\begin{example} \label{exa:poset}
    Let $\ddown = (1,3)$ and $\dup = (3,4)$.  Figure~\ref{fig:poset} gives both $\Pi_{\ddown,\dup}$ and its order complex. 
    Note that we use concatenation of the sequence to label the vertices, e.g. we label $(2,4)$ as $24$.
    \begin{figure}[!ht]
        \includegraphics{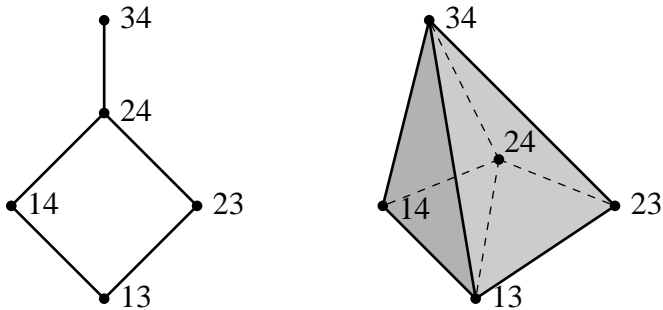}
        \caption{The Hasse diagram and the order complex of $\Pi_{(1,3),(3,4)}$}
        \label{fig:poset}
    \end{figure}
\end{example}

Further, the family of Boij-S\"oderberg posets contains the family of posets of bounded root sequences described in~\cite[Section~8]{ES}
where it is shown that the supernatural cohomology tables of root sequences in a bounded range give a geometric realisation of the order
complex of the associated bounded root sequences.  These are in turn used to prove the Boij-S\"oderberg conjectures.

\section{Basic structure} \label{sec:basic-structure}

Let $\ddown, \dup \in \ZZ_0^{p+1}$ be strictly increasing sequences with $\ddown_i \leq \dup_i$ for $0 \leq i \leq p$. First
we see that the order complex has very nice structure.
\begin{proposition}\label{pro:lattice}
    The poset $\Pi_{\ddown, \dup}$ is a bounded complete lattice.
\end{proposition}
\begin{proof}
    Let $D = \{d^{(1)}, \ldots, d^{(n)} \} \subset \Pi_{\ddown, \dup}$.  Define $l, u \in \ZZ^{p+1}$ by $l_i = \min_{j}d_i^{(j)}$ and
    $u_i = \max_{j}d_i^{(j)}$, for $0 \leq i \leq p$.  Then $l \leq d^{(j)} \leq u$ for $1 \leq j \leq n$.  Suppose $l' \leq d^{(j)}$
    for $1 \leq j \leq n$, then ${l'}_i \leq d_i^{(j)}$ for $0 \leq i \leq p$ and $1 \leq j \leq n$, so $l' \leq l$ and $l$ is the
    meet of $D$.  Similarly, $u$ is the join of $D$.
    
    Notice further that, by definition, $\ddown \leq d \leq \dup$ for all $d \in \Pi_{\ddown, \dup}$.  Thus, $\Pi_{\ddown, \dup}$
    is a bounded complete lattice.
\end{proof}

Next we see that reversing the order on $\Pi_{\ddown, \dup}$ yields a (possibly different) Boij-S\"oderberg poset.
\begin{proposition} \label{pro:reverse}
    Let $\widetilde{\Pi_{\ddown, \dup}}$ be the set $\Pi_{\ddown, \dup}$ endowed with the partial ordering given by
    $d \leq_{\sim} e$ if $d_i \geq e_i$ for $0 \leq i \leq p$.  Then $\widetilde{\Pi_{\ddown, \dup}}$ is isomorphic to some Boij-S\"oderberg poset.
\end{proposition}
\begin{proof}
    Without loss of generality, assume $\ddown_0 = 0$.  Let $m = \dup_p$ and define the map $\varphi$ from $\{0, \ldots, m\}$ onto itself by
    $i \mapsto m-i$ and extend $\varphi$ component-wise to $\{0, \ldots, m\}^{p+1}$.  Further, define the map $\rho$ from $\{0, \ldots, m\}^{p+1}$
    onto itself by $(d_0, \ldots, d_p) \mapsto (d_p, d_{p-1}, \ldots, d_0)$.  Notice $\rho\circ\varphi$ is a bijection.

    Then for $d, e \in \Pi_{\ddown, \dup}$, we have that
    \begin{eqnarray*}
        d \leq e & \Leftrightarrow & d_i \leq e_i, 0 \leq p \\
                 & \Leftrightarrow & \varphi(d_i) = m - d_i \geq \varphi(e_i) = m - e_i, 0 \leq p \\
                 & \Leftrightarrow & \rho(\varphi(d)) \geq \rho(\varphi(e)).
    \end{eqnarray*}
    That is, $\rho \circ \varphi$ exactly reverses the order of elements in $\Pi_{\ddown, \dup}$ and hence $\Pi_{\rho(\varphi(\dup)), \rho(\varphi(\ddown))}$
    is isomorphic to $\widetilde{\Pi_{\ddown, \dup}}$.
\end{proof}

An immediate result of this is that particular posets are isomorphic to themselves after reversing the order.
\begin{corollary}\label{cor:linear-symmetric}
    Let $b, m \in \NN$ and suppose $\ddown = (0, m, 2m, \ldots, pm)$ and $\dup = (b, m+b, 2m+b, \ldots, pm+b)$.  Then $\Pi_{\ddown, \dup}$
    is isomorphic to $\widetilde{\Pi_{\ddown, \dup}}$.
\end{corollary}
\begin{proof}
    Let $\varphi$ and $\rho$ be as in the proof of Proposition~\ref{pro:reverse}; recall that $\rho\circ\varphi$ is a bijection.
    Applying $\rho\circ\varphi$ to $\ddown$, we get
    \begin{eqnarray*}
        \rho(\varphi(\ddown)) &=& (pm + b - pm, pm + b - (p-1)m, \ldots, pm+b - m, pm+b - 0)\\
                              &=& (b, m+b, \ldots, (p-1)m+b, pm+b) \\
                              &=& \dup.
    \end{eqnarray*}
    Thus, $\rho\circ\varphi$ is a poset isomorphism, that is $\Pi_{\rho(\varphi(\dup)), \rho(\varphi(\ddown))} \cong \Pi_{\ddown, \dup}$.
\end{proof}

This then allows us to see that some posets can be simplified.
\begin{proposition}\label{pro:shearing}
    Suppose $\ddown = (0, \ldots, p)$ and $\dup = (k, \ldots, p+k)$ for some $1 \leq k \leq p$.  Then 
    \[
        \Pi_{\ddown, \dup} \cong \Pi_{(0, \ldots, k-1), (p+1, \ldots, p+k)}.
    \]
\end{proposition}
\begin{proof}
    By Corollary~\ref{cor:linear-symmetric}, $\Pi_{\ddown,\dup} \cong \widetilde{\Pi_{\ddown,\dup}}$.  Let
    $g: \binom{\{0, \ldots, p+k\}}{p+1} \rightarrow \binom{\{0, \ldots, p+k\}}{k}$ be given by $A \mapsto \{0, \ldots, p+k\} - A$.  
    
    Extending $g$ component-wise provides a bijective map from elements of  $\widetilde{\Pi_{\ddown, \dup}}$ to elements of 
    $\Pi_{(0, \ldots, k-1), (p+1, \ldots, p+k)}$.  Further still, for $d, e \in \widetilde{\Pi_{\ddown, \dup}}$, we have that $d \leq_{\sim} e$
    if and only if the $i^{th}$ largest missing element of $e$ is at least the $i^{th}$ largest missing element of $d$ for all
    $1 \leq i \leq k$.  That is, $d \leq_{\sim} e$ if and only if $g(d) \leq g(e)$.  Hence $g$ is a poset isomorphism.
\end{proof}

The binomial coefficients and the multi-dimensional Catalan numbers give the number of vertices and the number of facets, respectively,
for Boij-S\"oderberg posets with upper and lower sequences given by consecutive integers.
\begin{lemma} \label{lem:vertex-facet}
    Let $\ddown = (0, \ldots, p)$ and $\dup = (k, \ldots, p+k)$ for some positive integer $k$.  Then $\Pi_{\ddown, \dup}$ has
    $\binom{p+k+1}{p+1}$ vertices and
    \[
        f(p,k) := (pk+k)!\prod_{i=0}^{p}\frac{i!}{(k+i)!}
    \]
    facets.
\end{lemma}
\begin{proof}
    Every vertex of $\Pi_{\ddown, \dup}$ is a sequence of $p+1$ numbers which can be seen as a $p+1$ subset of $\{1, 2, \ldots, p+k\}$.  Similarly,
    every such $p+1$ subset can be seen as a strictly increasing sequence of $p+1$ numbers.  Hence there are $\binom{p+k+1}{p+1}$ vertices.

    Consider the vertices of $\Pi_{\ddown, \dup}$ with $\ddown$ subtracted from them, then we are considering weakly increasing sequences of length
    $p+1$ with entries from $\{0, 1, \ldots, k\}$.  Thus under the aforementioned consideration of the vertices, the maximal chains in $\Pi_{\ddown, \dup}$
    then correspond directly to the ${\rm SU}(p+1)$ walk diagrams of $(p+1)k$-steps as described in~\cite[Section~IV]{Di}.  Moreover, the number of
    such walk diagrams is given in~\cite[Equation~(4.8)]{Di} as $f(p,k)$.
\end{proof}

As every Boij-S\"oderberg poset contains and is contained in such a poset, we can give a bound on the number of vertices and number of facets
of an arbitrary Boij-S\"oderberg poset.
\begin{corollary}
    Let $\ddown, \dup \in \ZZ_0^{p+1}$ be strictly increasing sequences with $\ddown_i \leq \dup_i$ for $0 \leq i \leq p$.
    If we set $v$ to be the number of vertices and $n$ to be the number of facets of $\Pi_{\ddown, \dup}$, then
    \[
        \binom{\dup_0 - \ddown_p+2p+1}{p+1} \leq v \leq \binom{\dup_p-\ddown_0+1}{p}
    \]
    and
    \[
        f(p, \dup_0 - \ddown_p + p) \leq n \leq f(p, \dup_p - \ddown_0 - p).
    \]
\end{corollary}

We further make the observation that vertex-decomposability of $\De(\Pi_{\ddown, \dup})$ implies the vertex-decomposability
of any $\De(\Pi_{\ddown', \dup'})$ where $\ddown \leq \ddown' \leq \dup' \leq \dup$.
\begin{remark}\label{rem:vd-transfer}
    Let $d \in \Pi_{\ddown, \dup}$.  Then $\link_{\De(\Pi_{\ddown, \dup})}{d} = \Ga_d \star \Si_d$ where
    $\Ga_d = \De(\{e \in \Pi_{\ddown, \dup} \st e < d\})$ and $\Si_d = \De(\{e \in \Pi_{\ddown, \dup} \st d < e\})$,
    hence $d \star \Ga_d = \De(\Pi_{\ddown, d})$ and $d \star \Si_d = \De(\Pi_{d, \dup})$.

    For $\ddown \leq \ddown' \leq \dup' \leq \dup$, we then have that $\De(\Pi_{\ddown', \dup'})$ is obtained from
    $\De(\Pi_{\ddown, \dup})$ by linking, taking half of a join, and coning.  Thus, properties like vertex-decomposability
    and Cohen-Macaulayness, which respect linking, joining, and coning, are preserved.
\end{remark}

\section{Recursive structure} \label{sec:recursive-structure}

Boij-S\"oderberg posets admit recursive atom orderings which are simple to describe.
\begin{theorem}\label{thm:atom-ordering}
    Let $\ddown, \dup \in \ZZ_0^{p+1}$ be strictly increasing sequences with $\ddown_i \leq \dup_i$ for $0 \leq i \leq p$.
    Then $\Pi_{\ddown, \dup}$ admits a recursive atom ordering when the atoms are ordered lexicographically from smallest to largest.
\end{theorem}
\begin{proof}
    Let $\ddown, \dup \in \ZZ_0^{p+1}$ be strictly increasing sequences with $\ddown_i \leq \dup_i$ for $0 \leq i \leq p$; assume
    with out loss of generality that $\ddown_i < \dup_i$ for $0 \leq i \leq p$ (see~\cite[Lemma~3.3]{BS08a}).

    If $\ddown$ and $\dup$ differ in only one position, then $p = 0$ and clearly $\Pi_{\ddown, \dup}$ admits a recursive atom ordering.
    Suppose then that $\ddown$ and $\dup$ differ in more than one position, i.e., $p > 0$.  Define $e_i \in \ZZ^{p+1}$ by $(e_i)_j = 0$ 
    if $j \neq i$ and $(e_i)_i = 1$.
    
    The atoms of $\Pi_{\ddown, \dup}$ are exactly the elements $\ddown + e_i$ where $\ddown_i + 1 \leq \ddown_{i+1}$ or $i = p$; let
    $i_1 > \cdots > i_t$ be the indices of the atoms.  Then under the lexicographic ordering, $\ddown + e_{i_1} < \cdots < \ddown + e_{i_m}$.
    Furthermore, by induction on the number of positions where $\ddown$ and $\dup$ differ, we may assume that $P_j = \Pi_{\ddown + e_{i_j}, \dup}$
    admits a recursive atom ordering when the atoms are ordered lexicographically from smallest to largest.

    For $j \geq 2$, the atom $\ddown + e_{i_j} + e_{i_k}$ of $P_j$ covers $\ddown + e_{i_k}$ for $1 \leq k < j$ and these are the only atoms
    of $P_j$ which cover some $\ddown + e_{i_k}$ for $1 \leq k < j$.  Moreover, for $\ell > i_j$ not equal to some $i_k$ for $1 \leq k < j$,
    then $\ddown_\ell + 1 = \ddown_{\ell+1}$ so $\ddown + e_{i_j} + e_\ell$ is not a member of $P_j$.  Hence, the atoms $\ddown + e_{i_j} + e_{i_k}$,
    for $1 \leq k < j$, of $P_j$ are those which come first in the ordering of $P_j$.

    For $1 \leq k < j \leq t$, the atom $z = \ddown + e_{i_j} + e_{i_k}$ of $P_j$ is the join of $\ddown + e_{i_j}$ and $\ddown + e_{i_k}$ and
    covers both atoms.  Hence if $\ddown + e_{i_k}, \ddown + e_{i_j}<y$ for some $y \in \Pi_{\ddown, \dup}$, then $z \leq y$.

    Therefore $\Pi_{\ddown, \dup}$ admits a recursive atom ordering when the atoms are ordered lexicographically from smallest to largest.
\end{proof}

As the poset admits a recursive atom ordering, then the associated order complex is vertex-decomposable, hence squarefree glicci (so, in
particular, in the {\bf G}orenstein {\bf li}aison {\bf c}lass of a {\bf c}omplete {\bf i}ntersection, see~\cite{NR}) and Cohen-Macaulay.
\begin{corollary}\label{cor:vd}
    All Boij-S\"oderberg order complexes are vertex-decomposable, hence squarefree glicci and Cohen-Macaulay.
\end{corollary}
\begin{proof}
    Let $\ddown, \dup \in \ZZ_0^{p+1}$ be strictly increasing sequences with $\ddown_i \leq \dup_i$ for $0 \leq i \leq p$.  Then
    by Theorem~\ref{thm:atom-ordering}, $\Pi_{\ddown, \dup}$ admits a recursive atom ordering.  Thus by~\cite[Theorem~3.2]{BW83}
    and~\cite[Theorem~11.6]{BW97}, the order complex of $\Pi_{\ddown, \dup}$ is vertex-decomposable.

    Pure vertex-decomposable simplicial complexes are squarefree glicci~\cite[Theorem~3.3]{NR} and pure
    shellable~\cite[Theorem~2.8]{PB}, hence Cohen-Macaulay.
\end{proof}

However, not every Boij-S\"oderberg order complex remains in the family during shedding.
\begin{remark}
    Consider $\ddown = (1,3,4)$ and $\dup = (2,5,6)$ which is illustrated in Figure~\ref{fig:shedding}.  Then the shedding vertices of
    $\Pi_{\ddown,\dup}$ are $(1,3,6),$ $(1,4,5),$ $(1,5,6),$ $(2,3,4),$ $(2,3,6),$ and $(2,4,5)$.
    
    \begin{figure}[!ht]
        \includegraphics{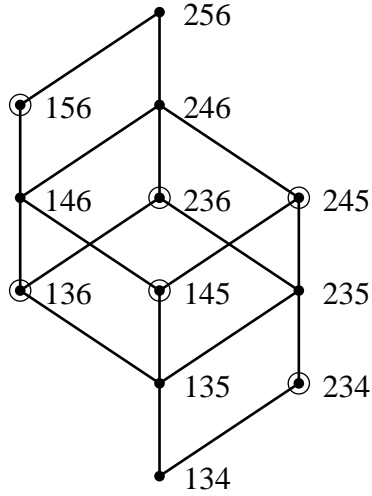}
        \caption{The Hasse Diagram of $\Pi_{(1,3,4),(2,5,6)}$ with shedding vertices circled}
        \label{fig:shedding}
    \end{figure}
    
    Removing either $(2,3,4)$ or $(1,5,6)$ creates a situation where a single change is followed by a triplet of changes
    or the opposite--either case is impossible in our family.  Last, removing the other four vertices creates a situation which is, after
    tedious calculation, demonstrable impossible.  In particular, the form implies $3 \leq \#d \leq 5$ but also that $\hat{0}$ would have
    exactly one consecutive pair of entries with nonconsecutive values.  Checking the nine possible situations yields a contradiction in
    each case.

    Hence, deleting any shedding vertex from $\Pi_{\ddown,\dup}$ yields a new poset which is not a Boij-S\"oderberg poset.
\end{remark}

\begin{acknowledgement}
    The author would like to thank his advisor, Uwe Nagel, for reading drafts of this article and making comments thereover.  The author would
    also like to thank Heide Gluesing-Luerssen for assistance in the aesthetics of the figures.
\end{acknowledgement}


\end{document}